\documentclass[12pt]{article}
\usepackage[cp1251]{inputenc}
\usepackage[english]{babel}
\usepackage{amssymb,amsfonts,amsmath}
\usepackage{amscd}
\newtheorem{Def}{Definition}[section]

\newtheorem{Th}{Theorem}[section]

\newtheorem{Lem}{Lemma}[section]
\newtheorem{Prop}{Proposition}[section]

\newenvironment{Proof}
{\par\noindent{\bf Proof.}}
{\hfill$\scriptstyle\blacksquare$}

\title{On infinite iterations of the functor of idempotent probability measures}
\author{Kh.~F.~Kholturaev\\
Tashkent institute of irrigation and agricultural mechanization engineers}

\begin{document}

\maketitle
\thispagestyle{empty}

\begin{center}

\end{center}
\begin{abstract}In this paper we establish that the functor of idempotent probability measures acting in the category of compacta and their continous mappings is perfect metrisable.\\


\end{abstract}

\tableofcontents

\section{\large Itroduction}

The notion of idempotent measure finds important applications in different parts of mathematics, mathematical physics, economics, mathematical biology and others. One can find a row of applications of idempotent mathematics from \cite{Litv2006arxiv}.

Let $\mathbb{R}$ be the real line. Consider the set $\mathbb{R}\cup \{-\infty\}$ with two algebraic operations: addition $\oplus$ and multiplication $\odot$ defined as $u \oplus v=\max\{u,\,v\}$ and $u\odot v=u+v$. The set $\mathbb{R}\cup \{-\infty\}$ forms semifield with respect to this operations and, the unity $\mathbf{1}=0$ and zero $\mathbf{0}=-\infty$, i.~e.
\begin{itemize}
\item[$(i)$] the addition $\oplus$ and the multiplication $\odot$ are associative;
\item[$(ii)$] the addition $\oplus$ is commutative;
\item[$(iii)$] the multiplication $\odot$ is distributive with respect to the addition $\oplus$;
\item[$(iv)$] each nonzero element $x\in \mathbb{R}\cup \{-\infty\}$ is intertible.
\end{itemize}
It denotes by $\mathbb{R}_{\max}$. It is idempotent, i.~e. $x\oplus x=x$ for all $x\in \mathbb{R}_{\max}$, and commutative, i.~e. the multiplication $\odot$ is commutative.

Let $X$ be a compact Hausdorff space, $C(X)$ the algebra of continuous functions $\varphi\colon X\to \mathbb{R}$ with the usual algebraic operations. On $C(X)$ the operations $\oplus$ and $\odot$ we define as following:
\begin{gather*}
\varphi \oplus \psi = \max\{\varphi,\,\psi\},\quad \mbox{where}\quad \varphi,\, \psi \in C(X),\\
\varphi \odot \psi = \varphi+\psi,\quad \mbox{where}\quad \varphi,\, \psi \in C(X),\\
\lambda \odot \varphi= \varphi + \lambda_X,\quad \mbox{where}\quad \varphi \in C(X),\quad \lambda \in \mathbb{R},
\end{gather*}
where $\lambda_X$ is a contant function.

Recall \cite{Zar2006arxiv30Aug} that a functional $\mu\colon\, C(X)\to \mathbb{R}$ is called to be an \textit{idempotent probability measure} on $X$, if:
\begin{itemize}
\item[$1)$] $\mu(\lambda_X)=\lambda$ for each $\lambda \in \mathbb{R}$;
\item[$2)$] $\mu(\lambda \odot \varphi)=\mu(\varphi)+\lambda$ for all $\lambda \in \mathbb{R}$, $\varphi \in C(X)$;
\item[$3)$] $\mu(\varphi \oplus \psi)=\mu(\varphi) \oplus \mu(\psi)$ for every $\varphi,\, \psi \in C(X)$.
\end{itemize}

For a compact Hausdorff space $X$ a set of all idempotent probability measures on $X$ we denote by $I(X)$. Consider $I(X)$ as a subspace of $\mathbb{R}^{C(X)}$. In the induced topology the sets of the view
\begin{gather*}
\langle \mu;\, \varphi_1,\, \dots,\, \varphi_k;\, \varepsilon \rangle=\{\nu \in I(X):\, |\mu(\varphi_i) - \nu(\varphi_i)| < \varepsilon, \, i = 1,\, \dots,\, k\},
\end{gather*}
form a base of neighbourhoods of the idempotent measure $\mu \in I(X)$, where $\varphi_i \in C(X)$,\, $i=1,\, \dots,\, k$, and $\varepsilon>0$. The topology generated by this base coincide with point-wise convergence topology on $I(X)$. The topological space $I(X)$ is compact \cite{Zar2006arxiv30Aug}. For a given map $f\colon\, X\to Y$ of compact Hausdorff spaces the map $I(f)\colon\, I(X)\to I(Y)$ defines by the formula $I(f)(\mu)(\varphi)=\mu(\varphi \circ f)$, $\mu\in I(X)$, where $\varphi \in C(Y)$. The construction $I$ is a covariant functor, acting in the category of compact Hausdorff spaces and their continuous maps.

Since $I$ is a normal functor, for an arbitrary idempotent measure $\mu \in I(X)$ we may define the support of $\mu:$\,\, $\mbox{supp}\, \mu = \cap \{A\subset X:\, \bar{A} = A,\, \mu \in I(A)\}$. For a point $x\in X$ by the rule $\delta_{x}(\varphi) = \varphi(x)$, $\varphi\in C(X)$, we  define the Dirac measure $\delta_{x}$ supported on the singleton $\{x\}$, i.~e. $\textrm{supp}\, \delta_{x} = \{x\}$.

Let $\mu_1,\, \mu_2 \in I(X)$.  Put
\begin{gather*}
\Lambda_{1\,2} = \Lambda (\mu_1,\, \mu_2) = \{\xi \in I(X^2):\, I(\pi_i)(\xi) = \mu_i,\, i= 1,\, 2\},
\end{gather*}
where $\pi_i \colon\, X \times  X \to X$ is the projection onto $i$-th factor, $i = 1,\, 2$. Note that $\Lambda_{1\,2} \ne \varnothing$ (see \cite{Zait2019arxiv}, Page 7).

Let $(X,\, \rho)$ be a metric compact space. A function $\rho_I\colon I(X)\times I(X) \to \mathbb{R}$ defined as
\begin{equation}\label{rhoI1}
\rho_I(\mu_{1},\, \mu_{2}) = \inf\{\sup\{\xi(\rho) \oplus \rho(x,\, y):\, (x,\, y)\in \mbox{supp}\,\xi\}:\, \xi\in \Lambda_{1\,2}\}.
\end{equation}
is a metric on $I(X)$ which is an extension of the metric $\rho$ (Theorem 1 \cite{Zait2019arxiv}). Besides, the metric $\rho_I$ generates point-wise converging topology on $I(X)$ (Theorem 2 \cite{Zait2019arxiv}). Note in \cite{Toj2010uzmj,ZaitToj2012arxiv15March} where considered another metric.

Note that a function
\begin{equation}\label{rhoI2}
\rho_I(\mu_{1},\, \mu_{2}) = \inf\{\sup\{\rho(x,\, y):\, (x,\, y)\in \mbox{supp}\,\xi\}:\, \xi\in \Lambda_{1\,2}\},
\end{equation}
suggested by A.~Zaitov, also is a metric on $I(X)$. It is obvious that the equalities (\ref{rhoI1}) and (\ref{rhoI2}) define the same metric.

\section{\large Preliminaries}

All of the concepts and results in this section we take from \cite{Fed1990Triples}. Under a functor we mean a covariant functor acting in the category $\mathfrak{Top}$ of topological spaces and their continuous maps, and some subcategories of $\mathfrak{Top}$.

\begin{Def}
{\rm A functor $\mathcal{F}$, acting in the category $\mathfrak{Comp}$ of compact Hausdorff spaces
and their continuous maps, is said to be \textit{seminormal} if it satisfies the following
conditions:
\begin{enumerate}
\item[$1)$] $\mathcal{F}$ preserves empty set and singleton, i.~e. $\mathcal{F}(\varnothing) = \varnothing$ and $\mathcal{F}(\mathbf{1}) = \mathbf{1}$ take places, where $\mathbf{1}$ is a singleton.
\item[$2)$] $\mathcal{F}$ preserves intersections, i.~e. $\mathcal{F}(\bigcap\limits_{F\in \mathcal{B}}F) = \bigcap\limits_{F\in \mathcal{B}}\mathcal{F}(F)$ for a given compact Hausdorff space $X$ and for every family $\mathcal{B}$ of closed subsets of $X$;
\item[$3)$] $\mathcal{F}$ is monomorphic, i.~e. $\mathcal{F}(i)\colon\, \mathcal{F}(A)\to \mathcal{F}(X)$ is an embedding for every given embedding $i\colon\, A\to X$;
\item[$4)$] $\mathcal{F}$ is continuous, i.~e. $\mathcal{F}(\lim \mathcal{S}) = \lim(\mathcal{F}(\mathcal{S}))$ for each spectrum $\mathcal{S} = \{X_{\alpha},\, \pi^{\beta}_{\alpha};\, A\}$ of compact Hausdorff spaces and their corresponding projections.
\end{enumerate}}
\end{Def}

If a functor $\mathcal{F}$ is seminormal then there exists a unique natural transformation $\eta^{\mathcal{F}} = \eta\colon\, \mbox{Id}\to \mathcal{F}$ of identity functor $\mbox{Id}$ into the functor $\mathcal{F}$. Moreover, this transformation is monomorphism, i.~e. for each compact Hausdorff space $X$ the map $\eta_{X}\colon\, X\to \mathcal{F}(X)$ is an embedding.

\begin{Def}
{\rm An acting in the category $\mathfrak{MComp}$ of metrisable compact spaces and their continuous maps seminormal functor $\mathcal{F}$ is said to be \textit{metrisable} if for any metrisable compact $X$ and for each metric $d = d_{X}$ on $X$ it is possible to assign a metric $d_{\mathcal{F}(X)}$ on the compact $\mathcal{F}(X)$ such that the following conditions hold:
\begin{enumerate}
\item[$P1)$] if $i\colon\, (X_{1},\, d^{1})\to (X_{2},\, d^{2})$ is an isometric embedding then
\begin{gather*}
\mathcal{F}(i)\colon\, (\mathcal{F}(X_{1}),\, d^{1}_{\mathcal{F}(X_{1})}) \to (\mathcal{F}(X_{2}),\, d^{2}_{\mathcal{F}(X_{2})})
\end{gather*}
also is an isometric embedding;
\item[$P2)$] the embedding $\eta_{X}\colon (X,\, d) \to (\mathcal{F}(X),\, d_{\mathcal{F}(X)})$ is an isometry;
\item[$P3)$] $\mbox{diam}(\mathcal{F}(X),\, d_{\mathcal{F}(X)}) = \mbox{diam}(X,\, d)$.
\end{enumerate}}
\end{Def}

Fix a seminormal functor $\mathcal{F}$ and a compact Hausdorff space $X$, and put $\eta_{n-1, n}=\eta_{\mathcal{F}^{n-1}(X)}:\mathcal{F}^{n-1}(X)\rightarrow \mathcal{F}^{n}(X)$. For positive integers $n < m$ put
\begin{gather*}
\eta_{n,m}=\eta_{m-1,m}\circ...\circ \eta_{n+1, n+2}\circ \eta_{n, n+1}.
\end{gather*}
A direct sequence
\begin{gather}\label{directsequen}
\begin{CD}
X @>{\eta_{0,\, 1}}>> \mathcal{F}(X) @>{\eta_{1,\, 2}}>> \dots @>{\eta_{n-1,\, n}}>>\mathcal{F}^{n}(X) @>{\eta_{n,\, n+1}}>> \mathcal{F}^{n+1}(X) @>{\eta_{n+1,\, n+2}}>> \dots \,\, .
\end{CD}
\end{gather}
arises.

Fix a metric $d$ on a compactum $X$ and a metrization of the functor $\mathcal{F}$. A metric on the compactum $\mathcal{F}^{n}(X)$, generated by this metrization we denote by $d_n$. Then every of the maps
\begin{gather*}
\eta_{n,m}:(\mathcal{F}^{n}(X), d_n)\rightarrow (\mathcal{F}^{m}(X), d_m)
\end{gather*}
is an isometric embedding. The limit of sequence (\ref{directsequen}) in the category metric spaces and their isometric maps we denote by $(\mathcal{F}^{+}(X),\, d_{+})$. We give more detail definition of the metric $d_{+}$.  Considering so far $\mathcal{F}^{+}(X)$ as a limit of (\ref{directsequen}) in the set category, a limit of the embeddings $\eta_{n,m}:\mathcal{F}^{n}(X)\rightarrow \mathcal{F}^{m}(X)$ as  $m\rightarrow \infty$ we denote by $\eta_n:\mathcal{F}^{n}(X)\rightarrow \mathcal{F}^{+}(X)$. Then
\begin{gather*}
\mathcal{F}^{+}(X)=\bigcup\{\eta_n(\mathcal{F}^{n}(X)): n\in \omega\},
\end{gather*}
and a metric $d_{+}$ defines by the metrics $d_{n}$ on the summands $\eta_{n}(\mathcal{F}^{n}(X))$, i.~e. for $x,\, y \in \mathcal{F}^{+}(X)$ we have
\begin{gather}\label{metr+}
\begin{CD}
d_{+}(x,\, y)=d_{n}(a,\, b)
\end{CD}
\end{gather}
where $\eta_{n}(a)=x$, $\eta_{n}(b)=y$. The definition of the metric $d_{+}$ by equality (\ref{metr+}) is correct, since at $n < m$ the maps $\eta_{n,\,m}$ are isometric embeddings.

If $f\colon\, X\to Y$ is a continuous map, then it is possible to determine a map $\mathcal{F}^{+}(f)\colon\, \mathcal{F}^{+}(X)\to \mathcal{F}^{+}(Y)$. It becomes as follows. For $x\in \mathcal{F}^{+}(X)$ there exist $n$ and $a\in \mathcal{F}^{n}(X)$ such that $x=\eta_{n}(a)$. We put $\mathcal{F}^{+}(f)(x) = \eta_{n}\mathcal{F}^{n}(f)(a)$. The correctness of this definition follows from what $\eta_{n,\, m}$ is a natural transformation of the functor $\mathcal{F}^{n}$ to the functor $\mathcal{F}^{m}$.

\begin{Def}\label{uniformdef}
{\rm A metrisable functor $\mathcal{F}$ is said to be \textit{uniformly metrisable}, if its some metrization has the following property
\begin{enumerate}
\item[$P4)$] for any continuous map $f\colon\, (X_{1},\, d^{1})\to (X_{2},\, d^{2})$ the map
\begin{gather*}
\mathcal{F}^{+}(f)\colon\, (\mathcal{F}^{+}(X_{1}),\, d^{1}_{+})\to (\mathcal{F}^{+}(X_{2}),\, d^{2}_{+})
\end{gather*}
is uniformly continuous.
\end{enumerate}}
\end{Def}

A metrization of a functor, satisfying property $P4)$, is called \textit{uniformly continuous}.

For a uniformly metrisable functor $\mathcal{F}$ the operation $\mathcal{F}^{+}$ is a functor, acting from the category $\mathfrak{MComp} = \mathfrak{Metr}\cap \mathfrak{Comp}$ of compacta into the category of (sigma-compact) metrisable spaces. Moreover Definition \ref{uniformdef} directly follows
\begin{Prop}
{\rm If $f \colon\, (X_{1},\, d^{1})\to (X_{2},\, d^{2})$ is a homeomorphism of compacta then for a uniformly metrisable functor $\mathcal{F}$ the map $\mathcal{F}^{+}(f)\colon\, (\mathcal{F}^{+}(X_{1}),\, d_{+}^{1}) \to (\mathcal{F}^{+}(X_{2}),\, d_{+}^{2})$ is a uniformly homeomorphism.}
\end{Prop}

Therefore a metric space $(\mathcal{F}^{+}(X),\, d_{+})$ topologically does not depend of the choose a metric $d$ on the compactum $X$. Consequently, the operation $\mathcal{F}^{+}$ may be consider as a functor, acting from the category compacta into the category metrisable spaces and uniformly continuous maps.

Each compactum $(X,\,d)$ is assigned a completion $(\mathcal{F}^{++}(X),\, d_{++})$ of the space $(\mathcal{F}^{+}(X),\, d_{+})$, and each continuous map $f\colon\, (X_{1},\, d^{1})\to (X_{1},\, d^{1})$ is assigned a map $\mathcal{F}^{++}(f)\colon\, (\mathcal{F}^{++}(X_{1})) \to (\mathcal{F}^{++}(X_{2}))$ which is an extension of the map $\mathcal{F}^{+}(f)$ on the completions of the spaces $(\mathcal{F}^{+}(X_{1}),\, d_{+}^{1})$ and $(\mathcal{F}^{+}(X_{2}),\, d_{+}^{2})$. Thereby it is defined a functor  $\mathcal{F}^{++}$, acting from the category metric compact spaces into complete metric spaces and uniformly continuous maps. The functor $\mathcal{F}^{++}$ is as topological invariant as the functor $\mathcal{F}^{+}$, and it is possible to consider $\mathcal{F}^{++}$ as a functor, acting from the category compacta into the Polish spaces.

Let for a seminormal functor $\mathcal{F}$ except the  natural transformation $\eta\colon\, \mbox{Id}\to \mathcal{F}$, it defines still a natural transformation $\psi^{\mathcal{F}}\colon\, \mathcal{F}^{2}\to \mathcal{F}$. If at the same time for every compact $X$ the  equalities
\begin{align}
\psi_{X}\circ \mathcal{F}(\eta_{X}) = \mbox{id}_{\mathcal{F}(X)},\label{compos1}\\
\psi_{X}\circ\eta_{\mathcal{F}(X)} = \mbox{id}_{\mathcal{F}(X)}, \nonumber
\end{align}
are executed, then the functor $\mathcal{F}$ is said to be \textit{semimonadic}. And if the equality
\begin{gather*}
\psi_{X}\circ \mathcal{F}(\psi_{X}) = \psi_{X}\circ \psi_{\mathcal{F}(X)}
\end{gather*}
is still carried out, then $\mathcal{F}$ is said to be \textit{monadic}.

For a positive integer $n$ we put
\begin{gather*}
\psi_{n+1,\, n} = \psi_{\mathcal{F}^{n-1}(X)}\colon\, \mathcal{F}^{n+1}(X)\to \mathcal{F}^{n}(X)
\end{gather*}
and for $n < m$:
\begin{gather*}
\psi_{m,\, n} = \psi_{n+1,\, n}\circ \psi_{n+2,\, n+1} \circ\dots \circ\psi_{m-1,\, m-2}\circ\psi_{m,\, m-1}.
\end{gather*}

For each semimonadic functor the following inverse sequence arises
\begin{gather}\label{invsequen}
\begin{CD}
\mathcal{F}(X) @<{\psi_{2,\, 1}}<< \mathcal{F}^{2}(X) @<{\psi_{3,\, 2}}<< \dots @<{\psi_{n,\, n-1}}<<\mathcal{F}^{n}(X) @<{\psi_{n+1,\, n}}<< \mathcal{F}^{n+1}(X) @<{\psi_{n+2,\, n+1}}<< \dots \,\, .
\end{CD}
\end{gather}
Denote by $\mathcal{F}^{\omega}$ the limit of the sequence (\ref{invsequen}). Since $\psi_{n+1,\, n}$ are natural transformations the operation $\mathcal{F}^{\omega}$ is functorial. The functor $\mathcal{F}^{\omega}$ acts as in the category $\mathfrak{Comp}$ as in its subcategory $\mathfrak{MComp}$.

For $1\le n < m$ put $q_{n,\, m} = \mathcal{F}(\eta_{n-1,\, m-1})$. Equality (\ref{compos1}) yields
\begin{gather*}
\psi_{m,\, n}\circ q_{n,\, m} = \mbox{id}_{\mathcal{F}^{n}(X)}.
\end{gather*}

\begin{Def}\label{perfectdef}
{\rm A uniformly metrisable semimonadic functor $\mathcal{F}$ is said to be \textit{perfect metrisable}, if some of its metrization along with properties $P1)\, -\, P4)$ has the following properties
\begin{enumerate}
\item[$P5)$] $\psi_{2,1}:(\mathcal{F}^{2}(X),\, d_{2}) \rightarrow (\mathcal{F}(X),\, d_{1})$ is a non-expending map;
\item[$P6)$] for every pair of $a\in \mathcal{F}^{2}(X)$ and $x\in X$ we have
\begin{gather*}
d_{1}(\psi_{2,1}(a),\, \eta_{0,1}(x)) = d_{2}(a,\, \eta_{0,2}(x)).
\end{gather*}
\end{enumerate}}
\end{Def}

Further, we denote by $\mathcal{F}$ a perfect metrisable functor. By $\theta_{n}\colon\, \mathcal{F}^{+}(X)\to \mathcal{F}^{n}(X)$ we denote a natural projection, which on the summand $\eta_{m}(\mathcal{F}^{m}(X))$ is defined as follows:
\begin{gather*}
\theta_n=\psi_{m,n}\circ \eta_{m}^{-1}.
\end{gather*}
The map $\theta_n$ non-expanding with respect to (\ref{metr+}) and $P5)$. Thence as uniformly continuous map $\theta_n$ extends on the completion $\mathcal{F}^{++}(X)$ of the space $\mathcal{F}^{+}(X)$. This extension we denote by $\overline{\theta}_{n}$.

At $n_{1} < n_{2}$ the following equality holds
\begin{gather}\label{tetan1}
\theta_{n_1}=\psi_{n_2,n_1}\circ \eta_{n_2}.
\end{gather}
Really, for $x\in \eta_{m}(\mathcal{F}^{m}(X))$ at $m>n_2$ we have
\begin{gather*}
\theta_{n_{1}}(x)=\psi_{m,\, n_1}\circ \eta_{m}^{-1}(x)=\psi_{n_{2},\, n_{1}}\circ \psi_{m,\, n_{2}}\circ \eta_{m}^{-1}(x)=\psi_{n_{2},\, n_{1}}\circ \theta_{n_{2}}(x).
\end{gather*}

Since a continuous map into a Hausdorff space is defined in a unique way with itself values on everywhere dense subspace, (\ref{tetan1}) implies
\begin{gather*}
\overline{\theta}_{n_1}=\psi_{n_2,n_1}\circ \overline{\theta}_{n_2}.
\end{gather*}
Therefore there exists unique map $\theta\colon\, \mathcal{F}^{++}(X)\to \mathcal{F}^{\omega}(X)$
(a limit of maps $\overline{\theta}_{n}\colon\, \mathcal{F}^{++}(X)\to \mathcal{F}^{n}(X)$) such that for each $n$ one has
\begin{gather*}
\overline{\theta}_{n}=\psi_{n}\circ \theta,
\end{gather*}
where $\psi_{n}\colon\, \mathcal{F}^{\omega}(X)\to \mathcal{F}^{n}(X)$ is a through projection of the inverse sequence (\ref{invsequen})

\begin{Prop}
{\rm Maps
\begin{gather*}
\eta_{n}\circ \overline{\psi}_{n}\colon\, \mathcal{F}^{++}(X)\to \mathcal{F}^{n}(X) \to \mathcal{F}^{+}(X) \hookrightarrow \mathcal{F}^{++}(X)
\end{gather*}
converge to the identity map uniformly on compact sets.
}
\end{Prop}

\begin{Prop}
{\rm For arbitrary pair of $x\in\mathcal{F}^{n}(X)$ and $a\in\mathcal{F}^{m}(X)$, $m\geq n+2$, we have
\begin{gather*}
d_{n+1}(\psi_{m,\, n+1}(a),\, \eta_{n,\, n+1}(x)) = d_{m}(a,\, \eta_{n, m}(x)).
\end{gather*}}
\end{Prop}

\begin{Th}
{\rm
For each compactum $X$ and every perfect metrisable functor $\mathcal{F}$ the map $\theta:\mathcal{F}^{++}(X)\rightarrow \mathcal{F}^{\omega}(X)$ is an embedding.}
\end{Th}

\begin{Th}
{\rm Let a perfect metrisable functor $\mathcal{F}$ and a compactum $X$ be such that
\begin{itemize}
\item[$1)$] $\mathcal{F}^{\omega}(X)$ is homeomorphic to $Q$;
\item[$2)$] $\mathcal{F}^{n}(X)$ is homeomorphic to $Q$ for all $n$, beginning with some;
\item[$3)$] $\eta_{n,\, n+1}(\mathcal{F}^{n}(X))$ is a $Z$-set in $\mathcal{F}^{n+1}(X))$ for all $n$.
\end{itemize}

Let, besides, the metrization of the functor $\mathcal{F}$ satisfies the following condition
\begin{enumerate}
\item[$P7)$] for an arbitrary $a\in \mathcal{F}(X)$ the inequality $d_1(a,\, \eta_{0,\, 1}(X))\ge \varepsilon$ implies the inequality $d_{i+1}(q_{1,\, i+1}(a),\, \eta_{i,\, i+1}q_{1,\, i}(a))\ge \varepsilon$ for all $i\ge 2$.
\end{enumerate}

Then the triple $(\mathcal{F}^{\omega}(X),\, \theta(\mathcal{F}^{++}(X)),\, \theta(\mathcal{F}^{+}(X)))$ is homeomorphic to the triple $(Q,\, s,\, \mbox{rint}\, Q)$}.
\end{Th}

Remind, that $Q = \prod\{[-1,\, 1]_{n}:\, n=1,\, 2,\, \dots\}$ is the Hilbert cube, $s$ is the pseudo-interior of $Q$, and $\mbox{rint}\, Q \, = \{x=\{x_n\}\in Q:\, |x_{n}|\le t < 1 \mbox{ for all } n=1,\, 2,\, \dots\}$. It is well-known that $s$ is homeomorphic to Hilbert space $l_{2}$, and $\mbox{rint}\, Q$ to $\Sigma$, which  is a linear span of $Q$ in $l_{2}$.

\section{\large Uniform metrizability of the functor of idempotent probability measures}

In this section we verify that the functor $I$ satisfies properties $P4) - P7)$. For this we need the following construction. Since functor $I$ is normal there exists unique natural transformation $\eta^{I} = \eta\colon\, \mbox{Id}\to I$ of identity functor $\mbox{Id}$ into the functor $I$. Here the natural transformation $\eta$ consists of monomorphisms $\delta_{X},\,X\in \mathfrak{Comp}$. More detail, the last means that for each compact Hausdorff space $X$ the map $\delta_{X}\colon\, X\to I(X)$, which defines as $\delta_{X}(x)=\delta_{x}$, $x\in X$, is an embedding. Thus $\eta = \{\delta_{X}:\, X\in \mathfrak{Comp}\}$ is the mentioned natural transformation.

Let $X$ be a compactum. Put
\begin{gather*}
I^{0}(X)=X,\quad I^{k}(X)=I(I^{k-1}(X)),\quad k=1,\, 2,\, \dots, \quad \mbox{ and }\\
\eta_{n-1,\, n} = \eta_{I^{n-1}(X)}\colon\, I^{n-1}(X)\to I^{n}(X).
\end{gather*}

Fix a metric $\rho$ on a compactum $X$ and consider the mertrization $\rho_{1}$ of the functor $I$ by the Zaitov metric $\rho_{1} = \rho_{I}$, defined by (\ref{rhoI2}). The metric on $I^{n}(X)$ generated by this metrization we denote as $\rho_{n}$.

\begin{Lem}\label{P1)}
{\rm Let $X$ be a compactum with a metric $\rho$. Then $\delta_{X}\colon\, (X,\, \rho)\to (I(X),\, \rho_{1})$ is an isometry.}
\end{Lem}

\begin{Proof}
For every pair of Dirac measures $\delta_{x} = 0\odot\delta_{x},\, \delta_{y} = 0\odot\delta_{y}\in I(X)$,  $x,\, y\in X$, the uniqueness of $(0\odot \delta_x,\, 0\odot\delta_y)$-admissible measure $0\odot \delta_{(x,\,y)}\in \Lambda(\delta_{x},\,\delta_{y})\subset I(X^{2})$ implies that
\begin{gather*}
\rho_{1}(\delta_x,\, \delta_y) = 0\odot \delta_{(x,\,y)}(\rho) \oplus \rho(x,\, y) = \rho(x,\,y).
\end{gather*}

\end{Proof}

\begin{Lem}\label{P2)}
{\rm Let $(X_{1}, \rho^{1})$, $(X_{2},\,\,\rho^{2})$ be compacta. Then for each isometric embedding $i\colon\, (X_{1}, \rho^{1}) \to (X_{2},\,\,\rho^{2})$ the map $I(i)\colon\,(I(X_{1}),\, \rho_{1}^{1})\to (I(X_{2}),\, \rho_{1}^{2})$ is also an isometric embedding.}
\end{Lem}

\begin{Proof}
Let $\mu_{1},\, \nu_{1}\in I(X_{1})$ and $\mu_{2} = I(i)(\mu_{1}),\,\nu_{2} = I(i)(\nu_{1})\in I(X_{2})$. Then $I(i)(X_{1}) \subset I(X_{2})$ and for every $\zeta\in \Lambda(\mu_{2},\, \nu_{2})$ supports on $(i\times i)\,(X_{1}\times X_{1})$. That is way we have
\begin{multline*}
\rho_{1}^{2}(\mu_{2},\, \nu_{2}) = \inf\{\sup\{\rho^{2}(x_{2},\, y_{2}):\, (x_{2},\, y_{2})\in \mbox{supp}\, \xi_{2}\}:\, \xi_{2}\in \Lambda(\mu_{2},\, \nu_{2})\} =\\
= \inf\{\sup\{\rho^{2}(i(x_{1}),\, i(y_{1})):\, (i(x_{1}),\, i(y_{1}))\in \mbox{supp}\, I(i\times i)(\xi_{1})\}:\, \xi_{1}\in \Lambda(\mu_{1},\, \nu_{1})\} = \\
= \inf\{\sup\{\rho^{1}(x_{1},\, y_{1}):\, (x_{1},\, y_{1})\in \mbox{supp}\,\xi_{1}\}:\, \xi_{1}\in \Lambda(\mu_{1},\, \nu_{1})\} = \rho_{1}^{1}(\mu_{1},\, \nu_{1}).
\end{multline*}

\end{Proof}

\begin{Lem}\label{P3)}
{\rm For any metric $\rho$ on the compactum $X$ the following equality holds
\begin{gather*}
\mbox{diam} (X,\,\rho) = \mbox{diam} (I(X),\,\rho_{1}).
\end{gather*}}
\end{Lem}

\begin{Proof}
The Proof immediately follows from (\ref{rhoI2}).

\end{Proof}

\begin{Lem}\label{munu-admiss}
{\rm \cite{Zait2019arxiv}. For every pair $\mu,\, \nu\in I(X)$ there exists $\lambda_{\mu\,\nu} \in \Lambda (\mu,\, \nu)$ such that $\rho_1(\mu,\, \nu) = \sup\{\rho(x,\, y):\, (x,\, y)\in \mbox{supp}\,\lambda_{\mu\,\nu}\}$.}
\end{Lem}

\begin{Lem}\label{calculdist}
{\rm Let $f\colon X\to Y$ be a continuous map of compacta $(X,\, \rho_{X})$, $(Y,\, \rho_{Y})$, and
\begin{gather*}
\rho_{X1}(\mu,\, \nu) = \sup\{\rho_{X}(x,\, y):\, (x,\, y)\in \mbox{supp}\,\lambda_{\mu\,\nu}\},
\end{gather*}
$\mu,\, \nu\in I(X)$. Then
\begin{gather*}
\rho_{Y1}(I(f)(\mu),\, I(f)(\nu)) \le \sup\{\rho_{Y}(f(x),\, f(y)):\, (x,\, y)\in \mbox{supp}\,\lambda_{\mu\,\nu}\}.
\end{gather*}}
\end{Lem}
\begin{Proof}
We put $\kappa = I(f\times f)(\lambda_{\mu\,\nu})$. We have
\begin{gather*}
I(q_{1})(\kappa)(\varphi) = I(q_{1})(I(f\times f)(\lambda_{\mu\,\nu}))(\varphi) = I(f\times f)(\lambda_{\mu\,\nu})(\varphi \circ q_{1}) = \\
= \lambda_{\mu\,\nu}(\varphi \circ q_{1}\circ f\times f) = \lambda_{\mu\,\nu}(\varphi \circ f\circ p_{1}) = \mu(\varphi \circ f) =I(f)(\mu)(\varphi),
\end{gather*}
i.~e. $I(q_{1})(\kappa) = I(f)(\mu)$. Here $p_{1}\colon X^{2}\to X$,  $q_{1}\colon Y^{2}\to Y$ are the projections onto the first factors, respectively, $\varphi\in C(Y)$. We have used the commutativity of the diagram
\begin{gather*}
\begin{CD}
X\times X @>{f\times f}>> Y\times Y\\
@V{p_{1}}VV @VV{q_{1}}V\\
X @>{f}>> Y ,
\end{CD}
\end{gather*}
i.~e. the equality $q_{1}\circ f\times f = f\circ p_{1}$. Similarly, one can show $I(q_{2})(\kappa) = I(f)(\nu)$, where $q_{2}\colon Y^{2}\to Y$ is the projection onto the second factor. Thus, $\kappa\in I(Y^{2})$ is a $(I(f)(\mu),\, I(f)(\nu))$-admissible measure. Hence, \begin{multline*}
\rho_{Y1}(I(f)(\mu),\, I(f)(\nu)) \le \sup\{\rho_{Y}(t,\, z):\, (t,\, z)\in \mbox{supp}\, \kappa\} = \\
=\sup\{\rho_{Y}(t,\, z):\, (t,\, z)\in \mbox{supp}\, I(f\times f)(\lambda_{\mu\,\nu})\} =\\
=\sup\{\rho_{Y}(f(x),\, f(y)):\, (x,\, y)\in \mbox{supp}\, \lambda_{\mu\,\nu}\}.
\end{multline*}

\end{Proof}

Note, that since the support $\mbox{supp}\,(^{1}\lambda_{^{1}\mu\,^{1}\nu})$ of $^{1}\lambda_{^{1}\mu\,^{1}\nu} = \lambda_{\mu\,\nu} \in I(X^2)$ is a compact subset of a given compact Hausdorff space $X$ one may write
\begin{gather*}
\rho_{X1}(^{1}\mu,\, ^{1}\nu) = \max\{\rho_{X}(x,\, y):\, (x,\, y)\in \mbox{supp}\,(^{1}\lambda_{^{1}\mu\,^{1}\nu})\}.
\end{gather*}
Therefore, there exist $(x_{^{1}\mu},\, y_{^{1}\nu})\in \mbox{supp}\,(^{1}\lambda_{^{1}\mu\,^{1}\nu})$ such that
\begin{gather}\label{calculbyrhox}
\rho_{X1}(^{1}\mu,\, ^{1}\nu) = \rho_{X}(x_{^{1}\mu},\, y_{^{1}\nu}).
\end{gather}

Let us consider $^{2}\mu$, $^{2}\nu$ $\in$ $I^{2}(X)$. Proposition \ref{calculdist} and equality (\ref{calculbyrhox}) consequently give:
\begin{gather}\label{1,0}
\rho_{X2}(^{2}\mu,\, ^{2}\nu) = \sup\{\rho_{X1}(^{1}\mu,\, ^{1}\nu):\, (^{1}\mu,\, ^{1}\nu)\in \mbox{supp}\,(^{2}\lambda_{\mu\,\nu})\},
\end{gather}
and
\begin{gather}\label{k2}
\rho_{X2}(^{2}\mu,\, ^{2}\nu) = \rho_{X1}(^{1}\mu_{^{2}\mu},\, ^{1}\nu_{^{2}\nu}) = \sup\{\rho_{X}(x,\, y):\, (x,\, y)\in \mbox{supp}\,(^{1}\lambda_{^{1}\mu\, ^{1}\nu})
\end{gather}
Similarly to the proof of inequality in Proposition \ref{calculdist}, one may establish
\begin{multline*}
\rho_{Y2}(I^{2}(f)(^{2}\mu),\, I^{2}(f)(^{2}\nu)) \le \\
\le \sup\{\rho_{Y1}(I^{1}(f)(^{1}\mu),\, I^{1}(f)(^{1}\nu)):\, (^{1}\mu,\, ^{1}\nu)\in \mbox{supp}\,(^{2}\lambda_{^{2}\mu\, ^{2}\nu})\}\\
\le \sup\{\rho_{Y}(f(x),\, f(y)):\, (x,\, y)\in \mbox{supp}\, (^{1}\lambda_{^{1}\mu\, ^{1}\nu})\},
\end{multline*}
i.~e.
\begin{gather}\label{k2f}
\rho_{Y2}(I^{2}(f)(^{2}\mu),\, I^{2}(f)(^{2}\nu)) \le \sup\{\rho_{Y}(f(x),\, f(y)):\, (x,\, y)\in \mbox{supp}\, (^{1}\lambda_{^{1}\mu\, ^{1}\nu})\}.
\end{gather}
We denote $\lambda_{2} = ^{1}\lambda_{^{1}\mu\, ^{1}\nu}$. Then (\ref{k2}) and (\ref{k2f}) has the following view
\begin{gather}\label{k2`}
\rho_{X2}(^{2}\mu,\, ^{2}\nu)  = \sup\{\rho_{X}(x,\, y):\, (x,\, y)\in \mbox{supp}\,\lambda_{2}\},
\end{gather}

\begin{gather}\label{k2f`}
\rho_{Y2}(I^{2}(f)(^{2}\mu),\, I^{2}(f)(^{2}\nu)) \le \sup\{\rho_{Y}(f(x),\, f(y)):\, (x,\, y)\in \mbox{supp}\, \lambda_{2}\}.
\end{gather}

\begin{Lem}\label{l_k}
{\rm Let $f\colon X\to Y$ be a continuous map of compacta $(X,\, \rho_{X})$, $(Y,\, \rho_{Y})$, and $^{k}\mu,\, ^{k}\nu\in I^{k}(X)$. Then for every $k\ge 1$ there exists a $\lambda_{k}\in I(X^{2})$ such that
\begin{gather*}
\rho_{X\,k}(^{k}\mu,\, ^{k}\nu) = \sup\{\rho_{X}(x,\, y):\, (x,\, y)\in \mbox{supp}\,\lambda_{k}\},
\end{gather*}
and
\begin{gather*}
\rho_{Y\, k}(I^{k}(f)(^{k}\mu),\, I^{k}(f)(^{k}\nu)) \le \sup\{\rho_{Y}(f(x),\, f(y)):\, (x,\, y)\in \mbox{supp}\,\lambda_{k}\}.
\end{gather*}}
\end{Lem}
\begin{Proof}
The proof immediately follows from  Lemma \ref{calculdist} and formulas (\ref{k2`} -- \ref{k2f`}) by induction.

\end{Proof}

\begin{Lem}\label{e-d}
{\rm If $f\colon\, X\to Y$ is a $(\varepsilon,\, \delta)$-uniformly continuous map of compacta, then for each $k\ge 1$ the map $I^{k}(f)\colon\, I^{k}(X)\to I^{k}(Y)$ is $(\varepsilon,\, \delta)$-uniformly continuous.}
\end{Lem}
\begin{Proof}
The proof immediately follows from the inequality in Lemma \ref{l_k}.

\end{Proof}

So, summarizing Lemmas \ref{P1)} -- \ref{e-d},  according to Definition \ref{uniformdef} we have the main result of the section.

\begin{Th}\label{uniformity}
{\rm Functor $I$ of idempotent probability measures is uniformly metrisable. More exactly, for the metrization, introduced by (\ref{rhoI2}), and for every continuous map $f\colon\, (X_{1},\, d^{1})\to (X_{2},\, d^{2})$ the map
\begin{gather*}
I^{+}(f)\colon\, (I^{+}(X_{1}),\, d^{1}_{+})\to (I^{+}(X_{2}),\, d^{2}_{+})
\end{gather*}
is uniformly continuous.}
\end{Th}

\section{\large Perfect metrizability of the functor of idempotent probability measures}

This section we begin the following statement.
\begin{Lem}
{\rm \cite{Zait2019arxiv} Every $\mu\in I(X)$ nay be represented as $\mu = \bigoplus_{x\in X} \lambda(x)\odot\delta_{x}\in I(X)$ in a unique way, where $\lambda\colon\, X\to \mathbb{R}$ is an upper semicontinuous map. For each $x_{0}\in X$ we have $x_{0}\in \mbox{supp}\, \mu$ if and of only if  $\lambda(x_{0}) > -\infty$.}
\end{Lem}

Consider now a system $\psi$. The system $\psi$ consists of all mappings $\psi_X\colon\, I^{2}(X)\to I(X)$, acting as the following. Given $M\in I^{2}(X)$ put $\psi_{X}(M)(\varphi) = M(\overline{\varphi})$, where for any function $\varphi \in C(X)$ the function $\overline{\varphi}\colon\, I(X)\to \mathbb{R}$ defines by the formula $\overline{\varphi}(\mu) = \mu(\varphi)$. Fix a compactum $X$ and for a positive integer $n$ put $\psi_{n+1,\, n}=\psi_{I^{n-1}(X)}\colon\, I^{n+1}(X)\to I^{n}(X)$. Note that $\psi_{n+1,\,n}\circ \eta_{n,\,n+1}=Id_{I^{n}(X)}$.

\begin{Lem}\label{P5)}
{\rm $\psi_{1,\,0}\colon\, (I^{2}(X),\, \rho_{2})\to (I(X),\, \rho_{1})$ is a non-expanding map.}
\end{Lem}
\begin{Proof} Let $M$, $N$ $\in$ $I^{2}(X)$. It is required to show that
\begin{gather}\label{rho1lerho2}
\rho_{1}(\psi_{1,\,0}(M),\, \psi_{1,\,0}(N)) \le \rho_{2}(M,\, N).
\end{gather}

By definition and equation (\ref{1,0}) there exists a $(M,\, N)$-admissible measure $\lambda_{MN}\in I(I(X)\times I(X))$ such that
\begin{gather*}
\rho_{2}(M,\, N) = \sup\{\rho_{1}(\mu,\, \nu):\, (\mu,\, \nu)\in \mbox{supp}\,\lambda_{MN}\},
\end{gather*}
For every $\varphi\in C(X)$ we have
\begin{gather*}
I(\pi_{I\,1})(\lambda_{MN})(\overline{\varphi}) = \lambda_{MN}(\overline{\varphi}\circ \pi_{I\,1}) = M(\overline{\varphi}) = \psi_{1,\,0}(M)(\varphi) = \overline{\varphi}(\psi_{1,\,0}(M)),
\end{gather*}
i.~e. $I(\pi_{I\,1})(\lambda_{MN})(\overline{\varphi}) = \overline{\varphi}(\psi_{1,\,0}(M))$. Analogously, $I(\pi_{I\,2})(\lambda_{MN})(\overline{\varphi}) = \overline{\varphi}(\psi_{1,\,0}(N))$. Here $\pi_{I\,i}\colon\, I(X)\times I(X)\to I(X)$ is the projection onto $i$-th factor. The last two equality imply $(\psi_{1,\,0}(M),\, \psi_{1,\,0}(N))\in \mbox{supp}\,\lambda_{MN}$. That is why
\begin{gather*}
\rho_{2}(M,\, N) = \sup\{\rho_{1}(\mu,\, \nu):\, (\mu,\, \nu)\in \mbox{supp}\,\lambda_{MN}\} \ge \rho_{1}(\psi_{1,\,0}(M),\, \psi_{1,\,0}(N)).
\end{gather*}

\end{Proof}

\begin{Lem}\label{P6)}
{\rm For each $N\in \psi^{-1}_{1,\,0}(\mu)$ we have $\rho_{1}(\mu,\,\delta_{x_0}) = \rho_{2}(\delta_{\delta_{x_0}},\,N)$.}
\end{Lem}
\begin{Proof}
The measure $\mu_{x_0}$ is a unique $(\mu,\, \delta_{x_0})$-admissible measure, where $\mbox{supp}\,\mu_{x_0} = \{(x,\, x_{0}):\, x\in  \mbox{supp}\,\mu\}$ and measures $\mu_{x_0}$ and $\mu$ act the same rule. So, we have
\begin{gather*}
\rho_{1}(\mu,\,\delta_{x_0}) = \sup\{\rho(x,\, x_{0}):\, (x,\, x_{0})\in \mbox{supp}\,\mu_{x_0}\} = \sup\{\rho(x,\, x_{0}):\, x\in \mbox{supp}\,\mu\}.
\end{gather*}

Similarly, the measure $N_{\delta_{x_0}}$ is a unique $(N,\, \delta_{\delta_{x_0}})$-admissible measure, and
\begin{gather*}
\rho_{2}(\delta_{\delta_{x_0}},\,N) = \sup\{\rho_{1}(\nu,\, \delta_{x_{0}}):\, \nu\in \mbox{supp}\,N\}.
\end{gather*}
Note that if $\psi_{1,\, 0}(N) = \mu$ then by definition $N(\overline{\varphi}) = \mu(\varphi)$ for all $\varphi\in C(X)$, from here $\mu\in \mbox{supp}\, N$. Therefore,
\begin{gather*}
\rho_{2}(\delta_{\delta_{x_0}},\,N) \ge  \rho_{1}(\mu,\, \delta_{x_{0}}).
\end{gather*}
Finally, inequality (\ref{rho1lerho2}) finishes the proof.

\end{Proof}

\begin{Lem}
{\rm If $\rho_{1}(\mu,\,\eta_{0,\,1}(X))\geq\varepsilon$ then $\rho_{2}(I(\eta_{0,\,1})(\mu),\,\,\eta_{1,\,2}(I(X)))\geq\varepsilon$.}
\end{Lem}
\begin{Proof}
It is clear that
\begin{gather*}
\rho_{1}(\mu,\,\eta_{0,\,1}(X)) = \rho_{1}\left(\mu,\,\bigoplus\limits_{x\in X} 0\odot\delta_{x}\right)
\end{gather*}
and
\begin{gather*}
\rho_{2}(I(\eta_{0,\,1})(\mu),\,\,\eta_{1,\,2}(I(X))) = \rho_{2}\left(\delta_{\mu},\,\bigoplus\limits_{\nu\in I(X)} 0\odot\delta_{\nu}\right).
\end{gather*}

For every $\varphi\in C(X)$ we have $\delta_{\mu}(\overline{\varphi}) = \overline{\varphi}(\mu) = \mu(\varphi)$ and
\begin{gather*}
\bigoplus\limits_{\nu\in I(X)} 0\odot\delta_{\nu}(\overline{\varphi}) =  \bigoplus\limits_{\nu\in I(X)}0\odot\nu(\varphi) = \bigoplus\limits_{x\in X} 0\odot\varphi(x) =  \bigoplus\limits_{x\in X} 0\odot\delta_{x}(\varphi).
\end{gather*}
Consequently, $\bigoplus\limits_{x\in X} 0\odot\delta_{x}\in \mbox{supp}\,\bigoplus\limits_{\nu\in I(X)}0\odot\nu$, which completes the proof.

\end{Proof}

Thus, we have established the following result.

\begin{Th}
{\rm The functor $I$ is perfect metrisable.}
\end{Th}

From above established results and results from works \cite{Fed1990Triples}, \cite{ZaitKhol2012uzmj} and  \cite{ZaitKhol2018arxiv} we have the main statement of the paper.

\begin{Th}
{\rm For every compact Hausdorff space $X$ containing more than one point the triple $(I^{\omega}(X),\, \theta(I^{++}(X)),\, \theta(I^{+}(X)))$ is homeomorphic to the triple $(Q,\, s,\, \mbox{rint}\, Q)$.}
\end{Th}

\section*{Acknowledgements}

The author would like to thank to professor Adilbek Zaitov – the head of the Department of Mathematics and Natural Disciplines of Tashkent Institute of Architecture and Civil Engineering for comprehensive support and attention.

\end{document}